\documentclass{icga}
\newcommand{\BibTeX}{{\rm B\kern-.05em{\sc i\kern-.025em b}\kern-.08em
    T\kern-.1667em\lower.7ex\hbox{E}\kern-.125emX}}

\usepackage{lscape}
\usepackage{rotating}

\usepackage{diago}
\def\substform#1#2{\hbox to 4.35\goheight{#1\hskip2pt\string @\hskip2pt#2\hss}}
\def\passform#1{\hbox to 3.5\goheight{#1\hskip2ptpass\hskip2pt\hss} }

\normalsize\golength11pt\gothinness0.4pt\gothickness1pt

\title{A CLASSIFICATION OF SEMEAI WITH APPROACH MOVES}
\runningtitle{A Classification of Semeai}

\author{Thomas Wolf\thanks{email:twolf@brocku.ca}}

\affiliation{Department of Mathematics and Statistics, 
Brock University, St Catharines, Ontario, Canada} 
\issue{December, 2015}

\begin{document}
\goboardini

\maketitle
\begin{abstract}
The paper analyzes liberty races (semeai) in the game of Go.
A rigorous treatment of positions with approach moves but no ko
revealed a class of positions not described in the literature
before. 

The complete presented classification is simpler and more compact 
than the one used so far. This was achieved through a classification 
by the types of liberties and the order they are to be filled instead 
of a classification by group status.

A principal difference to semeai without approach moves is that in the
presence of approach moves and shared liberties and the absence of eyes, 
two different strategies are possible because the un-orthodox filling of 
shared liberties instead of outer liberties may be a better strategy.

\end{abstract}

\section{Motivation} \label{motivation} 
In the game of Go, liberty races (semeai) belong to the few classes of
positions that can be analyzed without performing a tree search.
Roughly speaking, the reason is that the best moves of both sides do
not overlap and thus can be counted separately, apart from a possibly
coinciding capturing move in the shared region but that move can be
counted like other moves because it is the last move. 

Compared to numerous publications on the subject of semeai (see references) 
in this contribution we perform a complete classification of all
semeai that include ``plain'' approach moves (no fights, no kos) and
describe a case that has not been covered in other publications.

In \citebay{KFL82} Karl-Friedrich Lenz introduced the Semeai-formula
(see (\ref{KFLenz}) below) and described the basic principles of how
to evaluate semeai. The publications \citebay{RH96,RH98,RH03} of
Richard Hunter aim at Go players and contain many examples. The second
half of \citebay{RH03} discusses the appearance of ko in semeai but
not approach moves. In \citebay{TN2008,TN2009} Teigo Nakamura uses
Combinatorial Game Theory to add up the number of moves necessary to
occupy all outside liberties but does not consider the case that the
capturing move is outside {\em and} that approach moves are necessary.
The publications of Martin M\"{u}ller \citebay{MM99,MM01,MM03}
classify more general cases of semeai mainly from the point of view of
computer Go but with the consequence that exact relations can not be
formulated. The book \citebay{RJ1} of Robert Jasiek discusses basic
semeai including the case that one of the two essential chains is
inside an eye of the other essential chain but not the case of
approach moves. In \citebay{TR12} Thomas Redecker explains theoretical
foundations of a special type of Semeai which includes Ko.

As a by-product of our treatment of approach moves we are compacting
the classification of all cases to one table with one page of
instructions. Other classifications list many more cases. For example,
\citebay{RJ2} quotes Robert Jasiek about his book on semeai
``.. there's the five basic types of semeai with 93 possible cases and
over 200 principles governing how to determine status and outcome.''
(even without discussing approach moves). We think that a merger of
classes will not only be more satisfying for mathematically interested
readers and programmers but also for Go-players.

We will characterize the different cases and playing strategies by
equality relations, which we call 'balance relations', describing the
unsettled situation where strengths are balanced and success depends
on who moves first. If these equality-relations are not satisfied then the
inequality will directly tell who won unconditionally.

In order to minimize the number of cases and to reach a complete
classification in the presence of plain approach moves it is important not
to do the obvious thing and classify semeai phenomenologically by
possible outcomes (kill/kill or kill/seki) but by the source of their
different behavior which is the question where the capturing move is.
If a chain has an eye then there is no choice, the capturing move of
this chain has to be in the eye. But if the chain has no eye then the
capturing move can be in the inner region shared with the opponent
essential chain or it can be on one of the outside liberties as will
be shown by examples in section \ref{examples}.

The advantage of filling inside liberties first and having the capturing move 
outside is to be able to select that capturing move which avoids as many 
approach moves as possible (if there is such a choice).


The paper is organized as follows.
After defining basic semeai with approach moves in section~\ref{assume}
we introduce in section~\ref{relations} the notation and derive
equations of balance that characterize unsettled positions. 
Section~\ref{overview} gives an overview of all cases and
sections \ref{appmoves1}, \ref{appmoves2} discuss consequences of
approach moves.  Of special interest are semeai with a shared region 
and without eyes which are discussed in section~\ref{appmoves2}. 

The question when it is beneficial to convert a larger number of
approach moves into a smaller number of outside liberties is 
answered in section~\ref{choosing}. Section~\ref{OptimalPlay}
gives the proper order in which liberties are to be filled.
A generalization to multi-purpose approach moves is discussed in
section~\ref{general}. Section~\ref{summary} concludes 
the paper with a summary.

\section{Assumptions} \label{assume}
This paper considers positions in which only two chains (one white and
one black) are in a liberty race. Both are called essential
chains which are assumed to have several properties.
\begin{itemize}
\item None of the chains is under atari\footnote{A chain under atari 
  has only one liberty.} with the opponent moving next.
\item If one of the essential chains is captured then this settles the
  race.

  Examples where this is not the case are positions known as snapback,
  2-step snapback and 'under the stones' as described in the
  appendix.

  Another example for positions that are not considered are positions
  where an essential chain has a nakade shape, i.e.\ a shape that is
  small enough and compact enough so that the capture of that chain
  does not lead to two eyes for the other side (see {\tt
    http://senseis.xmp.net/?Nakade}).
\item Both essential chains are either in direct contact or separated
  by empty points that are neighbour to both chains.
\item Both chains can have up to one eye which must be small enough
  not to guarantee seki or life for the essential chain and it must be
  settled in the sense that it can not be split into two eyes in one
  move. Furthermore, the eye does not include an empty point that is
  not a liberty of the surrounding essential chain. Another assumption
  on the eye is that it does not contain all four points A1,A2,B1,B2
  plus another point in a corner. (If the opponent would occupy these
  four points and get captured then a move on B1 and A2 would give the
  opponent an eye in A1 which would change the move count in formula
  (\ref{eye}) below).

\item Both essential chains may require approach moves for their
  outside liberties. But these moves should be ``plain'' in the
  following sense. For example, if the essential chain is white then
  the position should not give White any gain or incentive to prevent
  Black from performing these approach moves, i.e.\ there is no ko and
  no other fight possible that could delay or preventing these moves.

  Furthermore, the sets of Black moves to approach different liberties
  of White should be disjoint. In other words there are no
  ``multi-purpose'' approach moves that serve to approach more than
  one white liberty.

  To summarize, it should be straight forward to count for each
  liberty of White how many approach moves by Black are required to
  occupy that liberty without the chance of the liberty filling move
  being back-captured.

  Generalizations to this strict definition of approach moves are
  discussed in sections~\ref{choosing},\ref{general}.

\end{itemize}
To summarize, all liberties (in the shared region, in the eyes and on
the outside) are ``plain'' and do not involve any fights and thus are
counted in integer numbers. 

These assumptions are chosen to be rather conservative in order to
allow the derivation of explicit formulas in the next section.
Minor generalizations are discussed in section \ref{general}.

\section{Relations characterizing unsettled Positions} \label{relations}
In this section we derive all relations that characterize unsettled
semeai positions satisfying the above criteria.  To be able to
describe all cases with only two mathematical relations we need to
introduce the following notation.

\subsection{Notation} \label{notation} 
We use Diagram~\ref{i2} to illustrate definitions that are made to describe
the position in Diagram~\ref{i1}.\footnote{Strictly speaking, the position 
in Diagram~\ref{i1} does not fully
satisfy the condition that approach moves to all liberties are
``plain''. White could play on g2 and have the strategy to
defend the new liberty f1 through ko. But because White has no eye and
has only one such liberty, our theory also applies to this position.
} Most of the variables introduced below
refer to properties of essential chains which are marked by
\twcrstone\ \tbcrstone\ in Diagram~\ref{i2}. If the variable has an
index then this specifies the colour of the chain. For example, if $Z$
is the size of an eye then $Z_B$ is the size of the black eye and
$Z_1$ is the size of the eye of player 1. (Player 1 will play the
role of what is called attacker in the literature in positions
without approach moves.)

Comments enclosed by [~~] in the following definitions 
refer to Diagram~\ref{i2}. 

\[
\myblock{15\goheight}{
\board1--17/1--6/
\alphafalse
\whites{D1E1K4L1L4M4N3N4O1O3Q1Q2P2P3O1F2F4G2G4H1H2H3H4J3}
\blacks{B1B2C2C3D3D4D5E2E5F1F3F5G5H5J4J5M1M2M3J2K1K2K3L3N2O2}
\shipoutgoboard}{$\!\!\!\!$}{\par \ \ }{i1}
\hspace{3\goheight}
\myblock{15\goheight}{
\board1--17/1--6/
\alphatrue
\whitecrosses{F2F4G2G4H1H2H3H4J3}
\blackcrosses{J2K1K2K3L3N2O2M1M2M3}
\whites{D1K4L1L4M4N3N4O1O3Q1Q2P2P3O1}
\blacks{B1B2C2C3D3D4D5E2E5F1F3F5G5H5J4J5}
\wletters{E1}
\wletters[9]{L1}
\letters[1]{C1}
\letters[1]{D2}
\letters[1]{E3}
\letters[1]{P1}
\letters[19]{J1}
\shipoutgoboard}{$\!\!\!\!$}{\par \ \ }{i2}
\]
\vspace*{-25pt}
\begin{tabbing}
$P_i$: \= \kill 
\>{\bf Input variables:}\\
$Z$: \= number \kill
$Z$: \> si\underline{\bf Z}e of an eye (number of points inside the eye 
        whether occupied or not) of an essential chain ($Z=0$ if no eye) \\
     \> [$Z_W=0$ (no eye), $Z_B=2$ (l1,l2)] \\
$I$: \> number of opponent stones \underline{\bf I}nside the eye 
        ($I=0$ if no eye) \\
     \> [$I_W=0$ (no eye), $I_B=1$ (stone \twlstone{I})] \\
$O$: \> number of external (\underline{\bf O}utside) liberties of an 
        essential chain \\
     \> [$O_W=3$ (e4, g3, g1), $O_B=1$ (n1)] \\
$A$: \> total number of \underline{\bf A}pproach moves to all outside 
        liberties of an essential chain, $A$ is an integer\\ 
     \> as we consider only ``plain'' approach moves [$A_W=4, A_B=1$ (points 
        marked A and \twlstone{A})] \\
$V$: \> number of approach moves to an outside liberty that can be 
        a\underline{\bf V}oided when the \\
     \> capturing move is on an outside liberty, i.e.\ 
        $V$ is zero if the chain has an eye, else $V$ is the\\
     \> maximal number of approach moves to any one of the outside liberties of 
        the essential chain\\
     \> [$V_W=3$ (c1, d2, e1), $V_B=0$ (because Black has an eye)] \\
$S$: \> number of \underline{\bf S}hared liberties [$S=1$ (point marked S)] \\
     \> \\
     \>{\bf Derived variables:} \\ 
${E}$: \> the difference of the number of moves played by both sides
        inside an \underline{\bf E}ye of size $Z$, with\\
     \> already $I$ opponent stones inside in order to capture the 
        chain with this eye \\
${X}$: \> the number of moves of the opponent to occupy safely all 
          e\underline{\bf X}clusive liberties of a chain, i.e.\\ 
       \> the number of all outside liberties and all their approach moves 
          and moves to fill the eye, if\\
       \>  necessary, repeatedly minus the $Z-3$ own moves inside the own 
          eye of size $Z$ if $Z>3$ to\\
       \> defend the eye, in other words $X$ is the number of moves to 
          capture a chain without the moves\\
       \> on shared liberties, under the assumption that the capturing 
          move is in the shared region 
        \\
${R}$: \> number of moves to capture a chain without the moves on shared space,
        all under the\\
       \> assumption that the capturing move is NOT in the shared region (a
        liberty count \underline{\bf R}educed by\\
       \> the number of approach moves to the outside capturing move)\\
$P_i$: \> \underline{\bf P}layers $P_1, P_2, P_I, P_{II} =$ Black/White 
\end{tabbing}
In the following we give relations between the symbols introduced above.  

\subsection{Derived Variables}
After determining input variables through a simple counting of liberties
and stones the other variables are determined as follows.  

The formula for the number ${E}$ of moves to fill an eye repeatedly
until it can be captured minus the number of own moves to defend the
eye is\footnote{We find it easier to remember for a given $Z$ to
  multiply $Z-1$ and $Z-2$ rather than the other form of the formula:
  $(Z^2-3Z+6)/2-I$ that appears more often in the literature. Go
  players may want to learn the values ${E}(Z)= Z \ \mbox{for}\ Z\leq
  3 \ \mbox{and}\ = 5,8,12,17 \ \mbox{for}\ Z=4..7$.}
\begin{equation} \label{eye}
 {E} = \left\{\begin{array}{ll} 
  {\rm for\ } Z=0,1 : & Z                    \vspace{2pt} \\
  {\rm for\ } Z>1   : & (Z-1)(Z-2)/2 + 2 - I .
            \end{array}
     \right.
\end{equation}

Applying  (\ref{eye})  to  Diagram~\ref{i2} we  get  ${E}_W=0$  because
$Z_W=0$ (White has no eye) and ${E}_B=(2-1)(2-2)/2+2-1=1$.

If a chain is captured and the capturing move is on an outside liberty
then one can  choose that outside liberty as  the capturing move which
requires the maximal number of approach moves because these approach 
moves need not be made.   The number $V$ of  approach moves that can 
be avoided this way is defined as
\begin{equation} \label{appmov}
 V = \left\{\begin{array}{ll}  {\rm for\ } Z\neq 0  : & 0 \vspace{6pt}
 \\ {\rm for\ } Z=0 : &  \mbox{max number of approach moves of any one
   of the outside liberties} 
            \end{array}
     \right.
\end{equation}
i.e.\ if the chain  has an eye then the capturing move  will be in the
eye, so $V=0$. For Diagram~\ref{i2} we get $V_W=3$ (the moves of Black
on c1,d2,e1  are needed  to prepare  Black on g1  whereas Black  on g3
needs only one  approach move on e3), and  $V_B=0$ because of $Z_B\neq
0$.

The formula for the number ${X}$ of all moves needed to occupy all
exclusive liberties of a chain (i.e.\ all liberties other than the
shared liberties) safely, i.e.\ including all approach moves if none 
of these moves is the capturing move, minus the own moves to defend the 
eye (if there is one) is
\begin{equation} \label{ExclusiveFill}
 {X} = {E} + O + A
\end{equation}
where ${E}$  is given by (\ref{eye}).  Applying (\ref{ExclusiveFill}) to
Diagram~\ref{i1} we get ${X}_W=0+3+4=7$ and ${X}_B=1+1+1=3$.

If one of the moves on an outside liberty is the capturing move then
$V$ moves less than counted under ${X}$ 
have to be played and the resulting number is
\begin{equation} \label{ExclusiveCapture}
 {R} = {X} - V
\end{equation}
with the values ${R}_W=7-3=4, \ {R}_B=3$ in Diagram~\ref{i2}.

The values ${X}, {R}$ and $S$ are used in the following to formulate
relations that characterize unsettled semeai positions.


\subsection{Relations of Balance} \label{BalEqn}
To  derive relations  that characterize  unsettled positions  we first
consider the case  that in an unsettled position  one such optimal and
successful sequence of  moves has its last move -  the capturing move -
in  the shared region.   Let us  assume White  was successful  in this
sequence.

Because White did  the capturing move in the  shared region, White has
to occupy all outer liberties  of Black and perform all approach moves
to  these liberties  before. That  means, up  to the  move  before the
capturing   move    White   had   to    do   ${E}_B+O_B+A_B+S-1={X}_B+S-1$
moves.\footnote{${E}_B$  is not needed  in this  formula because  if the
  capturing move is in the shared region then Black has no eye, but in
  that case ${E}_B=0$, so having ${E}_B$  in the formula is not wrong.} To
be successful and  do the capturing move, White had  to make the first
move in the sequence.  But because the position is unsettled and White
played  optimally, if  Black would  have started,  White should  not be
successful. The only way for White  not to be successful is that White
itself is  under atari one move  before White's capturing  move in the
shared region. Because  one empty point in the  shared region is left,
this is  a liberty of  White and it  is therefore the last  liberty of
White. That  means in  the last ${X}_B+S-1$  moves of Black,  Black must
have occupied all outside liberties of White and made all the approach
moves to  these liberties  but no moves  in the shared  region because
they  would have helped  White as  White's capturing  move was  in the
shared region. Because the  position is unsettled, the computed number
of moves of White and Black  must be equal so that whoever started and
can do  the next (capturing) move  wins. That means  ${X}_W\geq {X}_B$ and
${X}_W = {X}_B + S - 1$. In general the form of this relation of balance 
is\footnote{We call it relation and not equation because equations can
  be used to compute one variable in terms of all others and also
  because in our relations the $=$ sign may be replaced by $>$ or
  $<$.}
\begin{equation}
{X}_I = {X}_{II} + S - 1 \ \ \ \ \mbox{(FOF)} \label{FOF}
\end{equation}
where the inequality $({X}_I\geq {X}_{II})$ identifies players $I$ and $II$. 
We call this type FOF ({\bf F}ill {\bf O}utside liberties {\bf F}irst)
because before the capturing move all outside liberties are
filled.

The other case that none of the capturing moves of the two winning
sequences of Black and White is located in the shared region is easily
characterized. If an essential chain has no eye then optimal play will
select as capturing move that outer liberty which needs the most
approach moves because these approach moves need not be played if that
liberty is the capturing move.

Thus the second case is characterized by replacing in relation
(\ref{FOF}) the $S-1$ through $S$ because all $S$ shared liberties
are occupied before the capturing move and by replacing ${X}$ through
${X}-V={R}$ because not all approach moves are necessary if the chain to
be captured does not have an eye. Thus the relation of balance is 
\begin{equation}
{R}_1 = {R}_2 + S \ \ \ \ \mbox{(FIF)}  \label{FIF}
\end{equation}
where the inequality ${R}_1\geq {R}_2$ identifies 
players 1 and 2.\footnote{
In the (FIF) relation (\ref{FIF}) we use 1, 2 to label players
whereas in (\ref{FOF}) we label them with $I, II$ because 
$P_1$ need not be the same as $P_I$.}
The (FIF) relation (\ref{FIF}) applies also if player 1 or 2 or 
both have an eye because then $V=0$
for that player. We call this case FIF ({\bf F}ill {\bf I}nside
liberties {\bf F}irst)\footnote{More precise would be FILBCM 
({\bf F}ill {\bf I}nside {\bf L}iberties {\bf B}efore the {\bf C}apturing {\bf M}ove) but
  that would be hard to memorize.} because before the capturing move 
{\em all} inside liberties have to be filled by one player
which is player $P_1$.

To apply these relations to a position like Diagram~\ref{i2} we first
have to identify players I and II based on 
${X}_I\geq {X}_{II}$ before checking (\ref{FOF}) 
and players 1 and 2 based on ${R}_1\geq {R}_2$ 
before checking (\ref{FIF}). For Diagram~\ref{i2} we get 
$7={X}_W>{X}_B=3$ and thus 
\begin{equation} \label{be1}
7={X}_W>{X}_B+S-1=3
\end{equation}
for the FOF scenario and 
$4={R}_W>{R}_B=3$ and thus 
\begin{equation} \label{be2}
4={R}_W={R}_B+S=4 
\end{equation}
for the FIF approach. Relation (\ref{be1}) tells us as expected that
Black would be unsuccessful in occupying at first all exclusive
liberties of White by 7 moves, whereas 
relation (\ref{be2}) shows an equality, i.e.\
the position is unsettled, both sides can kill when moving first.

\[
\myblock{15\goheight}{
\board1--17/1--6/
\alphafalse
\blackfirstfalse
\whites{D1E1K4L1L4M4N3N4O1O3Q1Q2P2P3O1F2F4G2G4H1H2H3H4J3}
\blacks{B1B2C2C3D3D4D5E2E5F1F3F5G5H5J4J5M1M2M3J2K1K2K3L3N2O2}
\numbers{P1E4N1E3J1G3L2}
\shipoutgoboard}{\twpstone\ captures \\ \hspace*{3.5pt}}
{\par \ \ }{i3}
\hspace{3\goheight}
\myblock{15\goheight}{
\board1--17/1--6/
\alphafalse
\blackfirsttrue
\whitetris{E1}
\whites{D1K4L1L4M4N3N4O1O3Q1Q2P2P3O1F2F4G2G4H1H2H3H4J3}
\blacks{B1B2C2C3D3D4D5E2E5F1F3F5G5H5J4J5M1M2M3J2K1K2K3L3N2O2}
\numbers{E4P1E3N1G3}
\numbers[7]{C1}
\numbers[9]{D2}
\numbers[13]{G1}
\numbers[15]{J1}
\shipoutgoboard}{\tbpstone\ captures \\ \hspace*{3.6pt}
\tbstone{11}\ @ \twtstone,
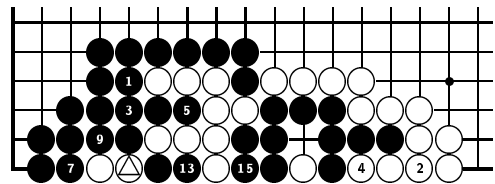
\twstone{6}, \twstone{8}, \twstone{10}, \twstone{12}, \twstone{14}\ 
pass \\ \hspace*{3pt}}
{\par \ \ }{i4}
\]

\subsection{Applicability of Relations of Balance} \label{AppBalEqn}
Because the first (FOF-)relation (\ref{FOF}) requires at least one
capturing move of a winning sequence in the shared region, this
relation is not applicable if there is no shared region and also not
if both sides have an eye, because then capturing moves can only be
moves in the eyes.

But even if only one side, say Black, has an eye
like in Diagram~\ref{i1}, the FOF-relation (\ref{FOF}) does not apply.
In that case a move capturing Black can only occur in Black's eye,
i.e.\ for White to have a chance it needs to occupy all shared
liberties.  If Black manages to occupy all but one outside liberty of
White (\tbstone{5} in Diagram~\ref{i4}) (in ${X}_W -V_W-1=7-3-1=3$
moves) before White occupies the last shared liberty (in ${X}_B+S-1 =
3+1-1=3$ moves which is equal to ${X}_B+S-1-V_B$ 
because $V_B=0$ due to the black eye) then White can occupy only
$S-1=1-1=0$ shared liberties before putting itself under atari
(i.e.\ without taking its last but one liberty and allowing to be
captured instantly). That means ${X}_W-V_W-1={X}_B+S-1-V_B$,
i.e.\ the (FIF-)relation (\ref{FIF}) applies and is the only one to 
apply if one side has an eye and the other not.


A comment about Diagram~\ref{i4}: The mathematical correctness of the
used formula is demonstrated in the above paragraph but one may wonder
why is the formula called FIF (Fill Inside First before capturing) if
the capturing move \tbstone{15}\ is inside and not outside. Answer: If
the only purpose is to win the liberty race then one of the best moves
for \twstone{6}\ is to play at \tbstone{15}. This move leads to
White's capture but all other moves would ultimately too.  At least
\twstone{6}\ on \tbstone{15}\ is needed for White to have any chance
to win. The FIF sequence is not the only best sequence of White,
another one also leading to capture is to pass, but FIF is one of the
best sequences and therefore it can be used to formulate a balance
relation and determine the outcome of the race. If White does not play
FIF (in order to save the move as a ko-threat for later) then Black is
allowed to change the reply too and play approach moves first before
capturing inside as shown in Diagram~\ref{i4}, or, even better, to
pass too and not to waste moves. But all of that does not change the
fact that purely for winning the liberty race, \twstone{6}\ on
\tbstone{15}\ is one of the best possible moves for White and
therefore suitable to establish a balance relation for determining the
status of the position.

In the remaining case of an existing shared region and no eyes both
relations of type FOF and FIF are relevant as discussed in section
\ref{caseB} further below.

The following section gives an overview of all cases and provides all
instructions to determine the status.

\section{Overview} \label{overview}
Table \ref{allcases} gives a complete overview of all cases of
semeai that involve approach moves and otherwise satisfy the 
assumptions of section \ref{assume}. We repeat both balance relations
(\ref{FOF}) and (\ref{FIF}) to have them close to the table:
\begin{eqnarray}
{X}_I & = & {X}_{II} + S - 1 \ \ \ \mbox{(FOF)} \label{FOFc} \\
{R}_1 & = & {R}_2 + S \ \ \ \ \ \ \ \ \ \ \ \ \, \mbox{(FIF)}.  \label{FIFc}
\end{eqnarray}
\vspace*{-12pt}
\begin{table}[!ht] 
\begin{center}
\begin{tabular}{|c|l|l|c|l|} \hline 
case &characterization of cases                           &def.\ of $P_i$& rel.& win of $P_1/P_2$ \\ \hline \hline
(A)&$S=0$                                                 & \ \ arbitrary        &(\ref{FIFc})& kill/kill \\ \hline
   &$S\neq 0,\  Z_W=Z_B=0$                                &\ \ (B)$_{\rm FOF}$:  &(\ref{FOFc})& kill/kill \ ($S=1$) \\ \cline{5-5}
(B)&                                                      &$\ \ {X}_I\geq{X}_{II}$&      & kill/seki ($S>1$)\\ \cline{3-5}
   &                                                      &\ \ (B)$_{\rm FIF}$:  &(\ref{FIFc})& kill/kill (played\\ 
   &                                                      &\ \ ${R}_1 \geq {R}_2$&       & only if $P_1$ can kill)\\ \hline
(C)&$S\neq 0,\ (0=Z_1<Z_2)$\ \ or                         &\ \ $Z_1<Z_2$         &(\ref{FIFc})& kill/kill \\ 
   &\ \ \ \ \ \ \ \ \ \ \,$(0<Z_1<Z_2>3)$                 &                      &       &           \\ \hline
(D)&$S\neq 0,\ \left((0<Z_1=Z_2)\right.$\ \ or            &\ \ ${X}_1 \geq {X}_2$&(\ref{FIFc})& kill/seki \\ 
   &\ \ \ \ \ \ \ \ \ \ \,$\left.(0<Z_1,Z_2\leq 3)\right)$&                      &       &           \\
\hline
\end{tabular}
\end{center}
\caption{Overview of all cases} \label{allcases}
\end{table}
\begin{center}
\end{center}

\newpage
\noindent{\bf Steps to determine the status using table \ref{allcases}:}
\begin{description} \itemsep=-1pt
\item[Step 1:] For a given position at first decide the case: 
  Case (A) has no shared liberties and case (B) has no eyes. In
  case (C) there is either only one eye or there are two eyes of
  different size of which the bigger one has a size $\geq
  4$.\footnote{In the literature eyes of size $\leq 3$ are called {\em
      small eyes} and eyes of size $\geq 4$ are called {\em big
      eyes}. The difference between both types is that for size $n\geq
    4$ capturing an inside opponent chain with $n-1$ stones costs one
    move which lets the opponent allow to make a move without cost but
    still $\geq 2$ liberties are left which is more than the one
    liberty before the capture. In other words, for size $\geq 4$
    capture buys time.} Case (D) has two eyes, either both of size
  $\leq 3$ or both of {\em equal} size $(\geq 4)$.
\item[Step 2:] Then identify players 1,2 based on inequalities 
    in column 2: \vspace{-4pt}
\begin{itemize}  \itemsep=0pt
  \item A distinction between players 1 and 2 is only needed for the
    allocation of shared liberties in the balance relation, i.e.\ in
    case (A) without shared liberties the distinction between player 1
    and 2 is irrelevant.
  \item In case (C) player 2 is identified as the one with an eye (if
    the other player has no eye) or, as the player with the bigger eye
    of size $\geq 4$. In case (D) player 1 has a higher or at least
    equally high sum of exclusive liberties plus approach moves 
    compared to player 2.
  \item In case (B) the player $I$ identified through ${X}_I\geq {X}_{II}$
    does not need to be the same as player 1 identified through 
    ${R}_1\geq {R}_2$ ! This case may need two steps. If player 1 
    identified through
    ${R}_1\geq {R}_2$ can kill based on relation (\ref{FIFc}) then this
    settles the status. Otherwise player $I$ has to be identified, now on
    the basis of ${X}_I\geq {X}_{II}$ with possibly a different result, and
    relation (\ref{FOFc}) is to be used to determine the status.  For
    example, in Diagram~\ref{p50} (further below) is $3={R}_B>{R}_W=1$ and
    therefore $P_1=$ Black whereas $4={X}_W>{X}_B=3$ would give $P_I =$
    White. Detailed explanations are given in section \ref{caseB}.
\end{itemize}
\item[Step 3:] The balance relations in the $3^{\rm rd}$ column 'rel.'
  describe unsettled positions for each case where the outcome of the
  race shown in the last column depends on who moves first. If the $=$
  is replaced by $>$ then player 1 wins always and if $=$ is replaced
  by $<$ then player 2 always wins. If the difference between both
  sides of the relation is $n$ then the loser is $n$ moves too slow.
\item[Step 4:] Optimal play is described in section \ref{OptimalPlay}.
\end{description}

\section{Comments} \label{comments}
The following comments concern the new case (B)$_{\rm FIF}$,
they relate known facts about semeai to our table  \ref{allcases}
and the relation of our terminology to 
definitions used in the literature.

\subsection{The new Case}
  From the two cases under (B) it is the case (B)$_{\rm FOF}$ that
  occurs most frequently in games. Case (B)$_{\rm FIF}$ has even not
  been described in the literature yet as far as the author knows. In
  section~\ref{caseB} it is shown that a necessary condition for being
  able to kill through FIF when reaching only seki under FOF is
  $V_2-V_1>1$ and when being killed under FOF is $V_2-V_1>2S+1$ where
  players 1,2 are identified through ${R}_1>{R}_2$. This shows that case
  (B)$_{\rm FIF}$ can only be relevant in the presence of approach
  moves and only if the winner under FOF has an outside liberty that
  needs at least 2 approach moves more than any outside liberty of the
  loser under FOF, or if one player can win through FIF and FOF but FIF allows
  more passes (see section \ref{bothkill}).

\subsection{Shared Liberties}
  In cases (B)$_{\rm FIF}$, (C) and (D) player 2 gets all $S$
  shared liberties in relation (\ref{FIFc}) and in case (B)$_{\rm FOF}$ 
  player 2 gets only $S-1$ shared liberties in the corresponding
  relation (\ref{FOFc}).

  In case (D) the player 2 with fewer exclusive liberties ${X}_2$
  obtains all shared liberties in the FIF formula because for player 2
  to be killed, player 1 needs to occupy all shared liberties,
  i.e.\ player 1 needs at least $S$ more exclusive liberties.


\subsection{Seki for Different Eye Sizes}
  If both sides have an eye of size $\leq 3$ then the status can
  still be seki even if both eyes are of different size. The reason is
  that for size 3 the catching of a 2-stone opponent chain inside the
  eye does not buy time because the opponent answers instantly and the 
  eye has still only one liberty inside. The situation is different when
  a 3-stone chain is caught inside an eye of size 4. Even if the opponent 
  puts a stone inside the eye immediately, the eye has at least 
  temporarily 2 inner liberties. This might be enough to occupy a shared 
  liberty and put the opponent under atari and win the semeai.

\subsection{Relation to other Publications}
\begin{itemize}
\item In table \ref{allcases} the positions with $S=1, Z_1=Z_2=0$ are
  listed under case (B) and not as usually in the literature under
  case (A). One could consider them under case (A) but case (A)
  contains eyes of all sizes whereas positions with $S=1, Z_1=Z_2=0$
  have no eyes. More importantly, these positions have the capturing
  move in the shared region and in the presence of approach moves they
  are characterized completely like case (B) positions. An example is
  shown in Diagrams \ref{p10} - \ref{p14} and discussed in section
  \ref{examples}.
\item In positions with shared liberties the literature on semeai
  identifies players as {\em attacker} and {\em defender}.
  In \citebay{MM99} and \citebay{RJ1} player 1 is called attacker and player
  2 is the defender, i.e.\ defender is a synonym for the side which
  gets the shared liberties in the liberty count.\footnote{Calling a
    player 'defender' who attacks as well, only from the outside, does
    not seem very obvious but that is what is used in the literature.}
  In \citebay{TN2008} on page 179, each external region has its defender
  (the owner of the essential chain) and its attacker (the other
  side), for example, in the external region of White the defender is
  White and Black the attacker. But then on page 180, the words
  attacker and defender are used like in all other publications.

  Players are also identified as {\em favourite} and {\em underdog}.
  In \citebay{RJ1},\citebay{RH96},\citebay{RH98},\citebay{RH03} the site with
  advantages (an eye versus no eye, or a bigger eye of size $\geq 4$
  versus a smaller eye, or in the case of no eyes then the side with
  more outer liberties) is called favourite and the other underdog.
  In kill/kill positions (C) player 2 is the favourite and
  gets the shared liberties due to its higher strength, in kill/seki
  positions (B)$_{\rm FOF}$, (D) the underdog gets the shared
  liberties because the underdog only wants to get a seki and does not
  try to kill.  Case (B)$_{\rm FIF}$ is not described in the
  literature.
  
  To summarize: What we call player~1 in cases (C),(D) is called 
  {\em  attacker} or {\em underdog} in other publications and what we call
  player~1 in case (B)$_{\rm FOF}$ is called {\em attacker} or {\em favourite}
  elsewhere.

  In this paper we do not identify players as attacker or defender
  because this would imply that one side can be identified as {\em
    the} attacker and the other as {\em the} defender which is not the
  case anymore in the presence of approach moves. For example, in
  Diagram~\ref{p52} further below Black plays the FOF approach and is
  the defender (i.e.\ shared liberties count for Black) whereas if
  Black plays the FIF approach in the same position as shown in
  Diagram~\ref{p53} then shared liberties count for White and thus
  Black would now have to be called attacker. Also the identification
  as favourite or underdog becomes problematic when in the absence of
  eyes one side has more exclusive physical liberties and the other
  requires more outside approach moves.
\item If one writes the balance equations in the form 
  \begin{eqnarray}
  \Delta_{\rm FOF}  & := & {X}_I - {X}_{II} = S - 1  \label{FOFd} \\
  \Delta_{\rm FIF}\ & := & {R}_1 - {R}_2 \ = S      \label{FIFd}
  \end{eqnarray}
  and defines what is called {\em forced liberties} in \citebay{MM99}
  \[F:=\left\{\begin{array}{l} S-1 \\
                               S   
              \end{array} \right. \ \ \  
              \begin{array}{l} {\rm for \ \ FOF} \\
                               {\rm for \ \ FIF}
              \end{array} \]
  then in the absence of approach moves 
  both relations (\ref{FOFd}), (\ref{FIFd}) become the Semeai-formula 
  \begin{equation} \label{KFLenz} \Delta = F 
  \end{equation}
  of Karl-Friedrich Lenz in \citebay{KFL82}, also shown 
  in \citebay{MM99} and \citebay{RJ1}.

\end{itemize}

\section{Approach Moves in Cases (A), (C) and (D)} \label{appmoves1}
As already explained in sections \ref{BalEqn} and \ref{AppBalEqn} the
generalization to approach moves is straightforward in cases (A), (C)
and (D). In these cases if a chain has no eye then the capturing move
is on an outside liberty and then an outside liberty can be chosen as
capturing move which otherwise would require the most approach
moves from all outside liberties. 

By defining the variable $V$ as the number of approach moves that can
be avoided (which is zero if the chain has an eye and otherwise is the
maximal number of approach moves to any outside liberty) and
introducing ${R}={X}-V$ we can use one and the same balance relation
(\ref{FIFc}) for all cases (A),(C),(D) whether chains have eyes or
not.

Because case (B) is more rich, a whole section is devoted to it.

\section{Approach Moves in Case (B)} \label{caseB}  \label{appmoves2}
In case (B) the situation is different as there are two different
types of play which both could be successful.  One play has the
capturing move inside, i.e.\ it {\bf F}ills all {\bf O}utside
liberties {\bf F}irst (FOF) - the standard play, and the other has the
capturing move outside and therefore {\bf F}ills {\bf I}nside
liberties {\bf F}irst (FIF) - the novel play.  

The following subsections discuss all aspects of case (B) and 
answer the following questions: \vspace{-3pt}
\begin{itemize} \itemsep=-3pt
\item What are necessary conditions for FIF to be better than FOF?
\item What are possible outcomes of both plays?
\item Can a win by FIF be enforced?
\item What if FIF and FOF can both kill?
\item Which steps decide whether to do FIF or FOF?
\item What are rules for optimal play? \vspace{-3pt}
\end{itemize}
We start with deriving necessary conditions for success with FIF if
FOF is sub-optimal.

\subsection{Necessary Conditions for FIF to be more useful than FOF}
This section is not necessary for any decision process but it
provides surprisingly short necessary conditions for a player to be 
more successful playing FIF than FOF. 

If FIF is played as described by formula (\ref{FIFc}) then player 1
has to occupy all inside liberties first. Player 1 would do this
only if he can win, otherwise he would play more defensively on
outside liberties first because in this case (B) both players have the
option to occupy outside liberties first. Therefore if we assume in
this subsection that for one player it is better to play FIF than 
FOF then this must be player 1 identified through ${R}_1>{R}_2$.

There are two cases:
Under FOF $P_1$ has seki (in the following case 1) or gets killed (case 2).

{\em Case 1:} \\
If we assume that $P_1$ moves first then we have
\begin{eqnarray}
{R}_1 & \geq & {R}_2+S   \label{nb1}\\
{X}_1 & <    & {X}_2+S-1 \label{nb2}
\end{eqnarray}
where (\ref{nb1}) follows from $P_1$ being successful under FIF and
(\ref{nb2}) follows from $P_1$ not being able to kill under FOF. 

This includes the two subcases that player 1 identified through
${R}_1>{R}_2$ is the same as either player $I$ or player $II$ identified through
${X}_I\geq {X}_{II}$. All that is assumed is that player 1 can not kill
through FOF which implies (\ref{nb2}).

Using at first (\ref{nb1}) together with ${R}_i={X}_i-V_i$ 
and then (\ref{nb2}) gives 
\begin{equation} \label{nb3} 
{X}_1-V_1={R}_1\geq {R}_2+S = {X}_2-V_2+S = ({X}_2+S-1)-V_2+1 > ({X}_1)-V_2+1
\end{equation}
and therefore
\begin{equation} \label{nb4}
 V_2-V_1 > 1 .
\end{equation}
If $P_2$ would move first then $\geq$ would become $>$
in (\ref{nb1}) and $<$ would become $\leq$ in (\ref{nb2}) and similarly in
(\ref{nb3}) but (\ref{nb4}) would not be affected.

For example, in Diagram~\ref{p20} below we have $P_1$ = White and
$V_1=0, V_2=1$, i.e.\ (\ref{nb4}) is not satisfied. White can kill
using FIF in Diagram~\ref{p24} {\em and} using FOF in Diagram
\ref{p22}, so FIF is not better.

Differently in Diagram~\ref{p50} where we have $P_1$ = Black and
$V_1=0, V_2=$, i.e.\ (\ref{nb4}) is satisfied. Indeed, Black can kill
using FIF in Diagram~\ref{p53} but reaches only a seki using FOF in
Diagram~\ref{p52}.

{\em Case 2:} \\
In the second case we assume that $P_1$ does even not reach seki under FOF
(i.e.\ $P_1=P_{II}$ because only $P_{II}$ can get killed under FOF) but 
$P_1$ can kill under FIF. If $P_1$ moves first then these two assumptions give
\begin{eqnarray}
{R}_1 & \geq & {R}_2+S \label{nb5}\\
{X}_2 & >    & {X}_1+S-1 \label{nb6}
\end{eqnarray}
and using at first (\ref{nb5}) and then (\ref{nb6}) we get
\begin{equation} \label{nb7}
{X}_1-V_1={R}_1\geq {R}_2+S=({X}_2)-V_2+S > ({X}_1+S-1)-V_2+S
\end{equation}
and therefore
\begin{equation} \label{nb8}
 V_2-V_1 > 2S-1 .
\end{equation}

This includes the two subcases that player 1 identified through
${R}_1>{R}_2$ is the same as either player $I$ or player $II$ identified through
${X}_I\geq {X}_{II}$. All that is assumed is that player 1 can not reach seki,
i.e.\ gets killed using FOF which implies (\ref{nb6}).)

If $P_2$ would move first then $\geq$ would become $>$
in (\ref{nb5}) and $>$ would become $\geq$ in (\ref{nb6}) and similarly in
(\ref{nb7}) but (\ref{nb8}) would not be affected.

For example, in Diagram~\ref{p70} we have $P_1$ = Black and $V_1=0,
V_2=2, S=1$, i.e.\ (\ref{nb8}) is satisfied. Indeed, Black can kill
playing FIF in Diagram~\ref{p73} even though it dies playing FOF in
Diagram~\ref{p72}.

To summarize: \\
It is possible to be able to kill with FIF when playing FOF does even not give 
seki (case 2) but it requires a bigger difference of the maximal number of
approach moves $V_2-V_1$ in (\ref{nb8}) than in (\ref{nb4}) for case 1 
when playing FOF at least gives seki. (In both inequalities $P_1$
is the winner of FIF.) In other words, playing FIF becomes attractive if 
the opponent has an outside liberty that requires many approach moves and
if oneself has only outside liberties that require no or only few approach 
moves. 

\subsection{Possible Outcomes of FIF and FOF and their Enforcement}
If all internal liberties are occupied before a capture happens then
there is no reason not to continue making moves; therefore FIF can not
result in seki.

Differently, when occupying outside liberties first, then unsettled
positions are of type kill/seki if $S>1$ (i.e.\ if there is more than
one inner liberty).

If player $P_1$ can kill through FIF 
then the other player can not prevent this. In other words, a FIF
kill by $P_1$ overwrites any FOF result, there is no competition in filling all
inside liberties first (FIF). Only one player can have an
interest in FIF which is player $P_1$ with ${R}_1>{R}_2$. If $P_2$ would
start filling all inner liberties then this would speed up the defeat
of $P_2$ and would be welcomed by $P_1$.

\subsection{Sensitivity of the decision between FIF and FOF}
As established above, only player $P_1$, i.e.\ the player with higher
$R$ value has the option to play FIF. Player $P_1$ will play FIF if
that wins the race. But whether $P_1$ wins may depend on who moves
first. For example, in Diagram~\ref{p20} below, player $P_1$ is White.
But White can win with FIF only if playing first, therefore if White 
plays first then White kills playing FIF as in Diagram~\ref{p24}. 
If Black plays first then White plays FOF as in Diagram~\ref{p21} 
and reaches seki. 

The decision whether player $P_1$ plays FIF or FOF may not only depend
on the position but also on the play of player $P_2$. 
For example, in Diagram~\ref{p80} player $P_1$ is Black
because $R_B=2>1=R_W$. If White, i.e.\ $P_2$ plays first then Black's
play depends on White's first move. If White decreases $X_B$ with
\twstone{1}\ in Diagram~\ref{p81} then Black should play FOF and if
White increases $V_W$ as in Diagram~\ref{p83} then Black should play
FIF.

\subsection{The Case that FIF and FOF can both kill}  \label{bothkill}
If player $P_1$ can kill through FIF and player $P_2$ can kill through
FOF\footnote{i.e.\ $P_I=P_2$ because only $P_I$ can kill in FOF, not $P_{II}$}
then $P_1$ can enforce the win as argued above. But what if $P_1$
can kill through FIF and FOF (i.e.\ $P_I=P_1$)?

If player $P_1$ can kill with both approaches then $P_1$
would want to pursue the more stable one, i.e.\ the approach where
$P_1$ can pass as often as possible.  Pursuing FIF, $P_1$ can pass
${R}_1-({R}_2+S)$ times according to (\ref{FIFc}) and pursuing FOF, $P_1$
can pass ${X}_1-({X}_2+S+1)$ times according to (\ref{FOFc}) (with $P_1=P_I$).  
Using ${R}_i={X}_i-V_i$ the difference between both passing numbers is
$D:=V_2-V_1-1$. That means if $D$ is positive then $P_1$ should FIF,
if $D$ is negative then $P_1$ should FOF and if $D=0$ then ko-fights
that are ongoing or may arise elsewhere on the board come into play.

If $P_1$ plays FOF then losing the semeai due to losing a ko-fight will still
get seki whereas pursuing FIF and losing the semeai due to ko would be a total
loss. On the other hand, it takes $P_1$ only ${R}_2+S$ moves to capture through
FIF but ${X}_2+S+1$ moves using FOF, i.e.\ pursuing FOF takes $V_2+1$ moves
more. This is the number of extra ko-threats that $P_2$ gets if $P_1$ plays
FOF. Thus, for $P_1$ to win a potential ko when playing FIF takes fewer but
larger ko-threats whereas winning a potential ko when playing FOF requires
more but lower value ko-threats.  That means the decision whether to play FIF
or FOF depends on the number of available ko-threats and their value in
comparison to the difference in killing and getting at least seki, and of
course it depends on whether $P_1$ can afford to play safe or not.

\subsection{Summary of Steps to choose between FIF and FOF} \label{decision}
For optimal play, both balance relations need to be checked. If one
player is successful in FIF this already decides the status, but for
practical play one still would have to check FOF as argued
above. Because success in killing through FIF is rare it is more
efficient to check at first FOF. The complete decision algorithm is
given through the following steps.
\begin{itemize}
\item Determine ${X}_B, {X}_W$ and player $P_I$ from ${X}_I\geq {X}_{II}$.
   Check the preliminary status based on FOF (\ref{FOFc}):
  if ${X}_I > {X}_{II}+S-1$ or ${X}_I = {X}_{II}+S-1$ and $P_I$ plays 
  first then $P_I$ doing FOF
  can kill otherwise if $S>1$ then it is a seki, if $S\leq1$ then
  $P_I$ gets killed. 
\item Determine $V_B, V_W$. If $V_B,V_W<2$ then stop, the FOF result is final.
\item Compute ${R}_B={X}_B-V_B,\ {R}_W={X}_W-V_W$ and determine player $P_1$
  from ${R}_1\geq {R}_2$. 
\item If ${R}_1<{R}_2+S$ or ${R}_1={R}_2+S$ and $P_2$ plays first then
  stop, the FOF result is final. (Independent of $P_1=P_I$ or
  $P_1=P_{II}$, it is for both players not advantageous to play
  FIF. Whoever plays inside liberties first gets killed.)
\item In the remaining cases $P_1$ kills. If $P_1$ can not kill if playing FOF 
  then it should play FIF.
\item In the remaining cases $P_1$ kills under FIF and FOF,
  i.e.\ $P_1=P_I$ (because only $P_I$ has a chance to kill under FOF).
  If possible, $P_1$ should play so to be able to pass and still win (if
  $V_2-V_1-1>0$ then FIF and if $V_2-V_1-1<0$ then FOF).
\item If $P_1$ can not pass {\em and} still win (i.e.\ $V_2-V_1-1=0$)
  then whether playing FIF or FOF does matter only for possible
  ko-fights. Playing FOF is less risky for $P_1$ (the outcomes are
  kill/seki instead of kill/kill) but gives $P_2$ more ko-threats
  because it takes $P_1$ more moves to catch $P_2$ under FOF.
\end{itemize}

\subsection{Examples}  \label{examples}
The following examples illustrate the difference between 
cases (B)$_{\rm FIF}$ and (B)$_{\rm FOF}$. Examples for
cases (A), (C), (D) are well known from the literature.
\[\hspace*{1\goheight}
\myblock{8.4\goheight}{
\board1--5/6--15/
\alphafalse
\whites{A7A8A12B11B12C7C8C12D8D9D10D11}
\blacks{A9A10A14B10B14C9C10C11C13C14D12D13}
\letters[1]{B8}
\shipoutgoboard}{$\!\!\!\!$}{\par \ \ }
{p10}
\myblock{8.4\goheight}{
\board1--5/6--15/
\alphafalse
\blackfirsttrue
\whites{A7A8A12B11B12C7C8C12D8D9D10D11}
\blacks{A9A10A14B10B14C9C10C11C13C14D12D13}
\numbers[1]{B13B8A13B9A11}
\shipoutgoboard}{$\!\!\!\!$}
{\par \ \tbpstone\ kills.}
{p11}
\myblock{8.4\goheight}{
\board1--5/6--15/
\alphafalse
\blackfirsttrue
\whites{A7A8A12B11B12C7C8C12D8D9D10D11}
\blacks{A9A10A14B10B14C9C10C11C13C14D12D13}
\numbers[1]{B13A11A13}
\shipoutgoboard}{$\!\!\!\!$}
{\par \ \tbpstone\ kills.}
{p12}
\myblock{8.4\goheight}{
\board1--5/6--15/
\alphafalse
\blackfirstfalse
\whites{A7A8A12B11B12C7C8C12D8D9D10D11}
\blacks{A9A10A14B10B14C9C10C11C13C14D12D13D14}
\numbers[1]{B8B13B9A13A11}
\shipoutgoboard}{$\!\!\!\!$}
{\par \ \twpstone\ kills.}
{p13}
\myblock{6.2\goheight}{
\board1--5/6--15/
\alphafalse
\blackfirstfalse
\whites{A7A8A12B11B12C7C8C12D8D9D10D11}
\blacks{A9A10A14B10B14C9C10C11C13C14D12D13}
\numbers[1]{A11B13B9}
\shipoutgoboard}{$\!\!\!\!$}
{\par \ \twpstone\ kills.}
{p14}\]
In Diagram~\ref{p10} we have $S=1$ and ${X}_W=2={X}_B$ and thus
balance condition (\ref{FOFc}) gives ${X}_W=2={X}_B+S-1$, thus the
position is unsettled, Black moving first, can kill as demonstrated in
Diagrams \ref{p11} and \ref{p12} and White can kill moving first as
shown in Diagram~\ref{p13}, both playing FOF. We also have 
${R}_W=2>1={R}_B$ and with FIF equation (\ref{FIFc})
further ${R}_W=2={R}_B+S$. Thus, White should 
also be successful when moving first and playing FIF by occupying 
all inner liberties before the outside capturing move which is shown 
in Diagram~\ref{p14}.
\[\hspace*{1\goheight}
\myblock{8.2\goheight}{
\board1--6/6--15/
\alphafalse
\whites{A7A8A12B12C7C8C12D8D9D10D11}
\blacks{A9A10A14B10B14C9C10C11C14D12D13D14}
\letters[1]{B8}
\shipoutgoboard}{$\!\!\!\!$}{\par \ \ }
{p20}
\myblock{8.2\goheight}{
\board1--6/6--15/
\alphafalse
\blackfirsttrue
\whites{A7A8A12B12C7C8C12D8D9D10D11}
\blacks{A9A10A14B10B14C9C10C11C14D12D13D14}
\numbers[1]{C13B8B13B9A13}
\shipoutgoboard}{$\!\!\!\!$}
{\par \ \tbpstone\ gets seki.}
{p21}
\myblock{8.2\goheight}{
\board1--6/6--15/
\alphafalse
\blackfirstfalse
\whites{A7A8A12B12C7C8C12D8D9D10D11}
\blacks{A9A10A14B10B14C9C10C11C14D12D13D14}
\numbers[1]{B8C13B9B13B11A13A11}
\shipoutgoboard}{$\!\!\!\!$}
{\par \ \twpstone\ kills.} 
{p22}
\myblock{8.2\goheight}{
\board1--6/6--15/
\alphafalse
\blackfirsttrue
\whites{A7A8A12B12C7C8C12D8D9D10D11}
\blacks{A9A10A14B10B14C9C10C11C14D12D13D14}
\numbers[1]{C13B11B13A11A13}
\shipoutgoboard}{$\!\!\!\!$}
{\par \ \tbpstone\ kills.}
{p23}
\myblock{7\goheight}{
\board1--6/6--15/
\alphafalse
\blackfirstfalse
\whites{A7A8A12B12C7C8C12D8D9D10D11}
\blacks{A9A10A14B10B14C9C10C11C14D12D13D14}
\numbers[1]{B11C13A11B13B9}
\shipoutgoboard}{$\!\!\!\!$}
{\par \ \twpstone\ kills.}
{p24}\]
Similarly in Diagram~\ref{p20}, now with $S=2$ and ${X}_W=3>2={X}_B$,
(\ref{FOFc}) gives: ${X}_W=2={X}_B+S-1$, i.e.\ the position is
unsettled when playing FOF: Black moving first gets seki 
(Diagram~\ref{p21}) and White moving first can kill (Diagram 
\ref{p22}). Here FIF gives ${R}_W=3>{R}_B=1$ and ${R}_W=3={R}_B+S$, 
i.e.\ also an unbalanced position (Diagrams~\ref{p23},~\ref{p24}) 
with the difference that the outcome for White in Diagram~\ref{p23} 
is worse than in Diagram~\ref{p21}, i.e.\ one should play FIF only 
if one can kill.
\[\myblock{9.4\goheight}{
\board1--6/5--17/
\alphafalse
\whites{A6A10A13A14B6B10B11B12C6C10D6D7D8D9}
\blacks{A8A12A16B8B13B15B16C8C9C11C12C13C14C15D10D11}
\letters[1]{A11}
\shipoutgoboard}{$\!\!\!\!$}
{\par \ \tbpstone\ to move.\\ $\ \ \ \ $}
{p50}
\myblock{9.4\goheight}{
\board1--6/5--17/
\alphafalse
\blackfirstfalse
\whitetris{A13}
\whites{A6A10A14B6B10B11B12C6C10D6D7D8D9}
\blacks{A8A12A16B8B13B15B16C8C9C11C12C13C14C15D10D11}
\numbers[1]{C7A15B7B14A7}
\numbers[7]{B9A11A9}
\shipoutgoboard}{\tbstone{6} @ \twtstone\\ }
{\par \twpstone\ kills.}
{p51}
\myblock{9.4\goheight}{
\board1--6/5--17/
\alphafalse
\blackfirsttrue
\whitetris{A13}
\whites{A6A10A14B6B10B11B12C6C10D6D7D8D9}
\blacks{A8A12A16B8B13B15B16C8C9C11C12C13C14C15D10D11}
\numbers[1]{A15C7B14B7}
\numbers[6]{A7A11}
\shipoutgoboard}{\tbstone{5} @ \twtstone\\ }
{\par \tbpstone\ gets seki.}
{p52}
\myblock{7\goheight}{
\board1--6/5--17/
\alphafalse
\blackfirsttrue
\whites{A6A10A13A14B6B10B11B12C6C10D6D7D8D9}
\blacks{A8A12A16B8B13B15B16C8C9C11C12C13C14C15D10D11}
\numbers[1]{B9C7A9B7A11}
\shipoutgoboard}{$\!\!\!\!$}
{\par \tbpstone\ kills.\\ $\ \ \ \ $}
{p53}\]
In Diagram~\ref{p50} the liberty A needs more than one approach move.
This is a necessary condition for FIF to be more useful than FOF for
the player that does not have the liberty A, here Black.  With $S=2$
and ${X}_W=4>3={X}_B$, FOF (\ref{FOFc}) gives: 
${X}_W=4={X}_B+S-1$, i.e.\ Black playing first can get seki 
(Diagram~\ref{p52}). The FIF approach gives ${R}_B=3>1={R}_W$ and with 
(\ref{FIFc}): ${R}_B=3={R}_W+S$, i.e.\ Black playing first can 
kill (Diagram~\ref{p53}). Therefore FIF is here better for Black.
\[\myblock{9.4\goheight}{
\board1--6/5--17/
\alphafalse
\whites{A6A10A13A14B6B10B11B12C6C10D6D7D9E9E8E7}
\blacks{A8A12A16B8B13B15B16C8C9C11C12C13C14C15D10D11}
\letters[1]{A11}
\shipoutgoboard}{$\!\!\!\!$}
{\par \ \tbpstone\ to move.\\ $\ \ \ \ $}
{p50a}
\myblock{9.4\goheight}{
\board1--6/5--17/
\alphafalse
\blackfirstfalse
\whitetris{A13}
\whites{A6A10A14B6B10B11B12C6C10D6D7D9E9E8E7}
\blacks{A8A12A16B8B13B15B16C8C9C11C12C13C14C15D10D11}
\numbers[1]{D8A15C7B14B7}
\numbers[7]{A7A11}
\shipoutgoboard}{\tbstone{6} @ \twtstone\\ }
{\par \twpstone\ gets seki.}
{p51a}
\myblock{9.4\goheight}{
\board1--6/5--17/
\alphafalse
\blackfirsttrue
\whitetris{A13}
\whites{A6A10A14B6B10B11B12C6C10D6D7D9E9E8E7}
\blacks{A8A12A16B8B13B15B16C8C9C11C12C13C14C15D10D11}
\numbers[1]{A15D8B14C7}
\numbers[6]{B7A11A7}
\shipoutgoboard}{\tbstone{5} @ \twtstone\\ }
{\par \tbpstone\ gets seki.}
{p52a}
\myblock{7\goheight}{
\board1--6/5--17/
\alphafalse
\blackfirstfalse
\whites{A6A10A13A14B6B10B11B12C6C10D6D7D9E9E8E7}
\blacks{A8A12A16B8B13B15B16C8C9C11C12C13C14C15D10D11}
\numbers[1]{D8B9C7A9B7A11}
\shipoutgoboard}{$\!\!\!\!$}
{\par \twpstone\ dies.\\ $\ \ \ \ $}
{p53a}\]
In Diagram~\ref{p50a} the black essential chain has one liberty more
than in Diagram~\ref{p50} leading to ${X}_W=4={X}_B$ and 
(\ref{FOFc}): ${X}_B=4<5={X}_W+S-1$,
i.e.\ to seki under FOF.  We also notice ${R}_B=4>1={R}_W$, i.e.\ the FIF
relation (\ref{FIFc}) gives ${R}_B=4>3={R}_W+S$. Therefore Black can kill
through FIF even if White plays first (Diagram~\ref{p53a}) whereas
under FOF even when moving first, Black can get only seki 
(Diagram~\ref{p52a}).

\[\myblock{9.4\goheight}{
\board1--6/4--17/
\alphafalse
\whites{B5A10A13A14B6B10B11B12C6C10D6D7D8D9}
\blacks{B7A8A12A16B8B13B15B16C8C9C11C12C13C14C15D10D11}
\letters[1]{A11}
\shipoutgoboard}{$\!\!\!\!$}
{\par \ \tbpstone\ to move.\\ $\ \ \ \ $}
{p60}
\myblock{9.4\goheight}{
\board1--6/4--17/
\alphafalse
\blackfirstfalse
\whitetris{A13}
\whites{B5A10A14B6B10B11B12C6C10D6D7D8D9}
\blacks{B7A8A12A16B8B13B15B16C8C9C11C12C13C14C15D10D11}
\numbers[1]{C7A15A6B14A7}
\numbers[7]{B9A11A9}
\shipoutgoboard}{\tbstone{6} @ \twtstone\\ }
{\par \twpstone\ kills.}
{p61}
\myblock{9.4\goheight}{
\board1--6/4--17/
\alphafalse
\blackfirsttrue
\whitetris{A13}
\whites{B5A10A14B6B10B11B12C6C10D6D7D8D9}
\blacks{B7A8A12A16B8B13B15B16C8C9C11C12C13C14C15D10D11}
\numbers[1]{A15C7B14A6}
\numbers[6]{A7A11}
\shipoutgoboard}{\tbstone{5} @ \twtstone\\ }
{\par \tbpstone\ gets seki.}
{p62}
\myblock{8\goheight}{
\board1--6/4--17/
\alphafalse
\blackfirsttrue
\whites{B5A10A13A14B6B10B11B12C6C10D6D7D8D9}
\blacks{B7A8A12A16B8B13B15B16C8C9C11C12C13C14C15D10D11}
\numbers[1]{B9C7A9A7}
\shipoutgoboard}{$\!\!\!\!$}
{\par \tbpstone\ dies.\\ $\ \ \ \ $}
{p63}\]
Diagram~\ref{p60} shows a position where both sides have approach
moves. Black has one outside liberty less and one approach move more,
so the same ${X}_B$ value. With $S=2$ and ${X}_W=4>3={X}_B$, FOF 
(\ref{FOFc}) gives
${X}_W=4={X}_B+S-1$, i.e.\ an unsettled position as demonstrated in
Diagrams \ref{p61}, \ref{p62}. The FIF approach gives ${R}_B=2>1={R}_W$
and with (\ref{FIFc}): ${R}_B=2<3={R}_W+S$ and therefore Black playing 
first dies under FIF in Diagram~\ref{p63}.

A comparison of Diagrams~\ref{p50} and \ref{p60} and their
corresponding results shows that a physical liberty and an approach
move have the same protection effect under FOF, i.e.\ both cost the
opponent a move. But for an attacking chain a physical liberty 
(of Black in Diagram~\ref{p50}) is more valuable than an extra
approach move (of Black in Diagram~\ref{p60}).
\[\myblock{8.4\goheight}{
\board1--5/5--16/
\alphafalse
\whites{A6A10A13A14B6B10B11B12C6C7C8C9}
\blacks{A8A12A15B8B9B13B15C10C11C12C13C14C15}
\shipoutgoboard}{$\!\!\!\!$}{\par \ \\ $\ \ \ \ $ }
{p70}
\myblock{8.4\goheight}{
\board1--5/5--16/
\alphafalse
\blackfirstfalse
\whites{A6A10A13A14B6B10B11B12C6C7C8C9}
\blacks{A8A12A15B8B9B13B15C10C11C12C13C14C15}
\numbers[1]{B7B14A7}
\numbers[5]{A9}
\shipoutgoboard}{$\!\!\!\!$}
{\par \ \twpstone\ kills.\\ $\ \ \ \ $}
{p71}
\myblock{8.4\goheight}{
\board1--5/5--16/
\alphafalse
\blackfirsttrue
\whitetris{A13}
\whites{A6A10A14B6B10B11B12C6C7C8C9}
\blacks{A8A12A15B8B9B13B15C10C11C12C13C14C15}
\numbers[1]{B14B7}
\numbers[4]{A7A11A9}
\shipoutgoboard}{\tbstone{3} @ \twtstone\\ }
{\par \ \tbpstone\ dies.}
{p72}
\myblock{6.4\goheight}{
\board1--5/5--16/
\alphafalse
\blackfirsttrue
\whites{A6A10A13A14B6B10B11B12C6C7C8C9}
\blacks{A8A12A15B8B9B13B15C10C11C12C13C14C15}
\numbers[1]{A9B7A11}
\shipoutgoboard}{$\!\!\!\!$}
{\par \ \tbpstone\ kills.\\ $\ \ \ \ $}
{p73}\]
Diagrams \ref{p70}-\ref{p73} give examples for the correctness of
relations  (\ref{FOFc}), (\ref{FIFc}) for $S=1$.
With $S=1,\ {X}_W=3>2={X}_2$ and FOF relation (\ref{FOFc}):
${X}_W=3>2={X}_B+S-1$
Black can not win even if moving first (Diagram~\ref{p72}).
Under FIF we get 
${R}_B=2>1={R}_W$ and with (\ref{FIFc}): ${R}_B=2={R}_W+S$, i.e.\ 
Black can kill if moving first (Diagram~\ref{p73}).
%

\[\hspace*{1\goheight}
\myblock{8\goheight}{
\board1--4/3--17/
\alphafalse
\blackfirstfalse
\whites{A4A10A14B4B5B10B11B12C5C6C7C8C9}
\blacks{A8A12A16B7B8B9B13B15B16C10C11C12C13C14C15}
\shipoutgoboard}{$\!\!\!\!$}{\par \ \ \par \ \ \par \ \ }
{p80}
\myblock{8\goheight}{
\board1--4/3--17/
\alphafalse
\blackfirstfalse
\whites{A4A10A14B4B5B10B11B12C5C6C7C8C9}
\blacks{A8A12A16B7B8B9B13B15B16C10C11C12C13C14C15}
\numbers{B6A13A6A11A7A9}
\shipoutgoboard}{$\!\!\!\!$}{\par Black plays\\FOF and kills\\}
{p81}
\myblock{8\goheight}{
\board1--4/3--17/
\alphafalse
\blackfirstfalse
\whites{A4A10A14B4B5B10B11B12C5C6C7C8C9}
\blacks{A8A12A16B7B8B9B13B15B16C10C11C12C13C14C15}
\numbers{B6A9A7}
\shipoutgoboard}{$\!\!\!\!$}{\par Black plays\\FIF and dies\\}
{p82}
\myblock{8\goheight}{
\board1--4/3--17/
\alphafalse
\blackfirstfalse
\whites{A4A10A14B4B5B10B11B12C5C6C7C8C9}
\blacks{A8A12A16B7B8B9B13B15B16C10C11C12C13C14C15}
\numbers{A13A9B6A11}
\shipoutgoboard}{$\!\!\!\!$}{\par Black plays\\FIF and kills\\}
{p83}
\myblock{6\goheight}{
\board1--4/3--17/
\alphafalse
\blackfirstfalse
\whites{A4A10A14B4B5B10B11B12C5C6C7C8C9}
\blacks{A8A12A16B7B8B9B13B15B16C10C11C12C13C14C15}
\numbers{A13B14B6A15A6A13A7A11A9}
\shipoutgoboard}{$\!\!\!\!$}{\par Black plays\\FOF and dies}
{p84}
\]
Diagrams \ref{p80} to \ref{p84} demonstrate that the decision by 
player $P_1$ whether to play FOF or FIF is taken when it is player $P_1$'s 
turn. $P_1$ is Black, because of $R_B=2>1=R_W$ in Diagram~\ref{p80}.
After \twstone{1}\ in Diagram~\ref{p81} is $X_B=2=X_W$ and $S=1$, i.e.\
(\ref{FOFc}) gives: $X_B=X_W+S-1$ therefore Black moving next can kill 
with FOF, whereas FIF would not work: (\ref{FIFc}) gives 
$R_B=1<2=1+1=R_W+S$ which is demonstrated in Diagram~\ref{p82}. 
If White plays \twstone{1}\ in Diagram~\ref{p83} then 
$V_W=3, X_W=4, R_W=1, V_B=1, X_B=3, R_B=2, S=1$, i.e.\ $R_B=R_W+S$ and 
Black playing next can kill with FIF as in Diagram~\ref{p83} whereas 
$X_B=3<4=X_W+S-1$, i.e.\ Black playing next dies with FOF as in 
Diagram~\ref{p84}.

\section{Choosing between Approach Moves and Outside Liberties}\label{choosing}
Diagram~\ref{ao1} shows only the outer region of a white essential
chain with some of its stones \twtstone. The question is whether the
shown part of the board already determines how many outer liberties
$O_W$, approach moves $A_W$ and maximal number $V_W$ of approach moves White
has.  If White would not move in its own outer region as in all other
positions of this paper then the values would be $O_W=1, A_W=5$ and
$V_w=5$ if White has no eye, else $V_W=0$.  As an alternative, White
could link its two chains like in Diagram~\ref{ao2} with the result
$O_W=2, A_W=V_W=0$. Which of both versions is better for White depends
on the hidden remainder of the position. In any case, the formulas
introduced in this paper can be used to decide which of both
alternatives is better for White.
\[
\myblock{9\goheight}{
\board1--8/1--3/
\alphafalse
\blackfirstfalse
\whitetris{f2g2h3h2h1}
\whites{b1c1d1e1}
\blacks{e2e3d3c3b3a3a2f1f3g3}
\shipoutgoboard}{$\!\!\!\!$}{\par \ \ }{ao1}
\ \ \ \ \ \ \ \ \ \ \ \ \ 
\myblock{9\goheight}{
\board1--8/1--3/
\alphafalse
\blackfirstfalse
\whitetris{f2g2h3h2h1}
\blackcrosses{f1}
\whites{b1c1d1e1}
\blacks{e2e3d3c3b3a3a2f3g3}
\numbers{g1a1}
\numbers[4]{b2}
\shipoutgoboard}{$\!\!\!\!$}{\par \twstone{3}\ at \tbcrstone.}{ao2}
\]

If White has an eye then $V_W=0$ and Diagram~\ref{ao1} is better for White,
because the capturing move would have to be in the eye, i.e.\ FIF 
is played requiring $O_W+A_W-V_W=6$ moves of Black to occupy all outside
liberties of White.

If White has no eye and either Black has an eye or there are no shared
liberties, then FIF is played as well but now $V_w=5$ in Diagram
\ref{ao1} with $A_W-V_W+Q_W=1$ move of Black is required to cover all
outside liberties of White.  In Diagram~\ref{ao2} White reaches
$A_W-V_W+Q_W=2$ which is slightly better in this situation.

Finally, if both sides have no eye (case (B)) in Table \ref{allcases}
and there is at least one shared liberty then FIF and FOF are
possible.  If FOF is played then Diagram~\ref{ao1} is (5+1)-(0+2)=4
moves better for White than Diagram~\ref{ao2} and if FIF is played
then Diagram~\ref{ao2} is 1 move better than Diagram~\ref{ao1}. FIF is
played only if one side filling inside liberties first can kill.

To decide whether to play as in Diagram~\ref{ao1} or as in \ref{ao2},
White has to check whether $O_W=1, A_W=5, V_w=5$ (Diagram~\ref{ao1})
would allow Black to kill with FIF (i.e.\ whether
$R_B(=O_B+A_B-V_B)>R_W+S(=O_W+A_W-V_W+S=1+S)$ or $R_B=R_W+S$ and Black
moves first) and then the only chance for White is to link as in
Diagram~\ref{ao2} and have one more liberty in FIF, i.e.\ hopefully to
have $R_B<R_W+S=2+S$ or $R_B=2+S$ and White with moving first, to
prevent FIF.

To summarize, if a player has the option to convert approach moves
into outside liberties then the Table \ref{allcases} with relations
(\ref{FOFc}) and (\ref{FIFc}) allow a quick decision about which 
variation is better for that player.

\section{Optimal Play} \label{OptimalPlay} Knowing the status, i.e.\
the optimal achievable result can already be useful: if the position
is settled and it is not urgent to eliminate ko-threats then one can
tenuki, i.e.\ play elsewhere and save moves as ko-threats if they are
at least forcing moves. On the other hand, if the position is
unsettled then one needs extra information about how to reach the
optimal result. The following are rules that are to be used in this
order for each move in unsettled positions.\footnote{Robert Jasiek is
  thanked for comments on exceptions for generalized semeai which 
  have been incorporated.}
These rules apply to all cases; only for case (B) in Table
\ref{allcases} one has to decide in criterion 5 based on section
\ref{decision} whether to do FIF or FOF.
\begin{itemize} \itemsep=3pt
\item[0.] Capture a chain of the opponent that is enclosing or approaching
  from the outside if that captures settles the semeai.\footnote{The
    reason for such a move to be more valuable than capturing the
    essential chain is that not only the semeai is won but more points
    and influence to the outside are gained. Strictly speaking, such
    an opportunity exists only for more general semeai with weak
    enclosing chains which are not treated in this paper. We still
    mention this type of move for the benefit of the Go-player.}
  Equivalently prevent the capture of an own chain that is enclosing
  or approaching from the outside if its capture would have settled
  the semeai.
\item[1.] Capture the opponent essential chain if it is under atari. 
\item[2.] If the own chain is under atari and has an eye which has an
  chain of the opponent inside that is under atari and has more than 2 stones
  then capture that opponent chain.
\item[3.] If the opponent just captured a chain in his eye allowing the
  opponent to get two eyes in the next move then prevent this
  eye splitting by playing on this point yourself.
\item[4.] Pass if the status is settled according to Table \ref{allcases}.
\item[5.] Fill or approach an outside liberty of the opponent chain. 
  Exception to
  this rule only for case (B)$_{\rm FIF}$: Do not occupy the last
  outside liberty and its approach moves which is the liberty with the
  most approach moves from all opponent's outside liberties.
\item[6.] Make a move in the eye of the opponent in a way that your 
  throw-in stones have a dead eye shape (nakade shape). 
\item[7.] Make a move on a shared liberty.
\end{itemize}

\section{Generalizations} \label{general}
The following generalizations may change the way how moves are counted
but not the balance relations (\ref{FOFc}), (\ref{FIFc}) themselves.

The assumption that each approach move is associated with only one
outside liberty was only made to be able to determine $V$ as it is
defined in section \ref{notation} easily as the maximum number of
approach moves to any one outside liberty of the essential chain. If
approach moves are not independent or serve to approach different
outside liberties then $V$ needs to be determined not through
(\ref{appmov}) but directly as number of approach moves to an outside
liberty that can be avoided if the capturing move is on an outside
liberty.  All other statements about $V$ remain unchanged, i.e.\ its
definition and its appearance in all relations.

\section{Summary} \label{summary}
In this paper we provide a complete analysis of semeai positions under
the assumptions of section \ref{assume}) which apart from the common
restrictions for semeai allow an arbitrary number of approach moves
without ko.

We derive relations of balance, i.e.\ mathematical conditions which
characterize positions where strengths are balanced, in other
words, positions that are unsettled where each side moving first can be
successful. If positions are settled then these relations give the
number of moves the losing side is behind which can be useful when
fighting kos.

We distinguish between the case that the capturing move is in the
shared region (when both players {\bf F}ill {\bf O}utside liberties
{\bf F}irst - called 'FOF') and the case that the capturing move is
not in the shared region (when one player has to {\bf F}ill {\bf I}nside
liberties {\bf F}irst - called 'FIF', i.e.\ before the capturing move)
and derive relations for both cases. By classifying positions
according to these two cases we manage to fit the whole classification
into only one table.

In this approach we also show that the case of no eyes and a single
shared liberty should be taken as a sub-case of the case of shared
liberties in order to generalize naturally the treatment to approach
moves whereas in the literature the main distinction is made between
positions with status kill/kill and those with status kill/seki which
merges the case of no shared region with 1-point shared regions.

By performing a complete derivation and not only presenting an ad hoc
collection of positions we find a case that has not been analyzed in
the literature yet: two chains with a shared region and both without
eyes.  For these situations we are able to derive very simple
necessary conditions for the case that the new FIF approach is
successful when the traditional FOF approach gives only seki or even
death.

A result, that at least for the author was unexpected, is that for
semeai without eyes but with approach moves and with a shared region
(even when consisting of only one point), one can not avoid to check
two cases.  The check is simple and even simpler necessary conditions
are available that help, but nevertheless, two cases need to be
considered. This is in contrast to basic semeai without approach moves
that can always be decided in a single decision process.  In this case
of shared liberties and no eyes only the player $P_1$ with more
external liberties plus approach moves minus the highest number of
approach moves to one external liberty has the option to play FIF.
$P_1$ will play FIF only if $P_1$ wins and that may depend on whether
$P_1$ moves first and if the opponent $P_2$ moves first then
whether $P_1$ can win by playing FIF may depend on the move $P_2$ 
plays first.

All statements are illustrated in numerous examples.

There are several generalizations that are of interest:
\begin{itemize}
\item The generalization to more than two essential chains.
\item The merger of the work in \citebay{TN2008,TN2009} with our work on
  approach moves.
\item The inclusion of ko as a type of approach moves.
\end{itemize}

\section*{Acknowledgment}
The author thanks Robert Jasiek, Martin M\"{u}ller, Harry Fearnley
and the referee for comments to the final manuscript.

\bibliography{semeai_icga}

\appendix
Typical for all the positions listed in this appendix is that the
white essential chain (on the right of each diagram) although being
able to capture the essential black chain on its left, can still not
link to the living enclosing white chain on the left.

\section{Snapback} \label{mouse1}
\[
\myblock{6\goheight}{
\board16--19/1--3/
\alphafalse
\whites{Q1Q2Q3R1S2T2}
\blacks{R2R3S3T3}
\shipoutgoboard}{$\!\!\!\!\!\!\!\!\!\!$}{}{mouse11}
\hspace*{2.5\goheight}
\myblock{7\goheight}{
\board16--19/1--3/
\alphafalse
\blackfirstfalse
\whites{Q1Q2Q3R1S2T2}
\blacks{R2R3S3T3}
\blacktris{S1}
\shipoutgoboard}{$\!\!\!\!\!\!\!\!\!\!$}{}{mouse12}
\hspace*{2.5\goheight}
\myblock{7\goheight}{
\board16--19/1--3/
\alphafalse
\blackfirstfalse
\whites{Q1Q2Q3R1S2T2}
\blacks{S1R2R3S3T3}
\numbers{T1}
\shipoutgoboard}{$\!\!\!\!\!\!\!\!\!\!$}{}{mouse13}
\hspace*{2.5\goheight}
\myblock{7\goheight}{
\board16--19/1--3/
\alphafalse
\blackfirstfalse
\whites{Q1Q2Q3R1}
\blacks{R2R3S3T3}
\numbers[2]{S1}
\shipoutgoboard}{$\!\!\!\!\!\!\!\!\!\!$}{}{mouse14}
\] 
Positions known as 'snapback' like the one in Diagram~\ref{mouse11}
are not semeai positions because there is no black essential
chain. But after the first move \tbtstone\ in Diagram~\ref{mouse12}
the position is at least topologically a semeai.  It is listed in this
appendix because White has no chance despite catching the black
essential chain.

\section{2-Step Snapback} \label{mouse2}
\parbox{0.57\textwidth}{
\[
\myblock{7\goheight}{
\board2--8/1--3/
\alphafalse
\whites{B1B2B3C1D2E2F2G2G1}
\blacks{D1E1C2C3D3E3F3G3H1H2H3}
\shipoutgoboard}{$\!\!\!\!\!$}{}{mouse21}
\hspace*{1.1\goheight}
\myblock{7\goheight}{
\board2--8/1--3/
\alphafalse
\blackfirstfalse
\whites{B1B2B3C1D2E2F2G2G1}
\blacks{D1E1C2C3D3E3F3G3H1H2H3}
\numbers{F1}
\shipoutgoboard}{$\!\!\!\!\!$}{}{mouse22}
\hspace*{1.1\goheight}
\myblock{7\goheight}{
\board2--8/1--3/
\alphafalse
\blackfirstfalse
\whites{B1B2B3C1D2E2F1F2G2G1}
\blacks{C2C3D3E3F3G3H1H2H3}
\numbers[2]{D1}
\shipoutgoboard}{$\!\!\!\!\!$}{}{mouse23}
\] 
\[
\myblock{7\goheight}{
\board2--8/1--5/
\alphafalse
\whites{B3B4B5C3D4E4F4G4G3B2B1D2E2F2G2}
\blacks{D3E3C4C5D5E5F5G5H3H4H5C2C1D1E1F1G1H2H1}
\shipoutgoboard}{$\!\!\!\!\!$}{}{mouse24}
\hspace*{1.1\goheight}
\myblock{7\goheight}{
\board2--8/1--5/
\alphafalse
\blackfirstfalse
\whites{B3B4B5C3D4E4F4G4G3B2B1D2E2F2G2}
\blacks{D3E3C4C5D5E5F5G5H3H4H5C2C1D1E1F1G1H2H1}
\numbers{F3}
\shipoutgoboard}{$\!\!\!\!\!$}{}{mouse25}
\hspace*{1.1\goheight}
\myblock{7\goheight}{
\board2--8/1--5/
\alphafalse
\blackfirstfalse
\whites{B3B4B5C3D4E4F4G4G3B2B1D2E2F2G2F3}
\blacks{C4C5D5E5F5G5H3H4H5C2C1D1E1F1G1H2H1}
\numbers[2]{D3}
\shipoutgoboard}{$\!\!\!\!\!$}{}{mouse26}
\]
\[\myblock{6\goheight}{
\board1--9/1--6/
\alphafalse
\blackfirstfalse
\whites{A3A4A5B3B5C3C5D3D5E2E4F2F4G2G4H2H3H4}
\blacks{A2B2C2D2D1E1F1G1H1J1J2J3J4J5H5G5F5E5E6D6C6B6A6C4D4E3F3}
\shipoutgoboard}{$\!\!\!\!\!$}{}{mouse27}
\hspace*{1.1\goheight}
\myblock{6\goheight}{
\board14--19/1--6/
\alphafalse
\blackfirstfalse
\whites{P1Q1P2P3R2S2T2}
\blacks{R1S1Q2Q3R3S3T3}
\shipoutgoboard}{$\!\!\!\!\!$}{}{mouse28}
\hspace*{1.1\goheight}
\myblock{6\goheight}{
\board14--19/1--6/
\alphafalse
\blackfirstfalse
\whites{S2S3T3R1Q1Q2Q3}
\blacks{S1T1R2R3R4S4T4}
\shipoutgoboard}{$\!\!\!\!\!$}{}{mouse29}
\]
}
\hspace*{0.04\textwidth}
\parbox{0.38\textwidth}{In these positions the black essential chain can
  afford to be captured twice and Black is still able to capture the white
  capturing chain afterwards. The white chain can not be linked to its
  strong white chain on the left.  This is possible because the black
  captured chain includes only 2 stones initially. From the first
  capture a 'normal' snapback results.

  The position in Diagram~\ref{mouse29} is special in that when White 
  captures the two black stones at t2, Black plays
  atari at s1, but then the position is "torazu san moku", such that
  a) if Black captures, White can recapture the 2 black stones,
  and either get an eye or capture another black stone, or b)
  if White captures, then a possible ko results. (Diagrams 
  \ref{mouse27} to \ref{mouse29} had been provided by 
  Harry Fearnley.)}

\section{Under the Stones} \label{under}
The technical term 'Under the Stones' refers to a situation where 
after the first capture the cutting stone \tbstone{2}\ is not captured
but can capture White in the next move. A possible characterization is:
\begin{itemize}
\item After capture, the capturing chain has only 2 liberties both
  resulting from the capture (if the capturing chain  would have afterwards
  only one liberty then it would just be a single ko or a snapback).
\item If after the capture the captured cross-cutting stone is replaced
  (by \tbstone{2}\ in diagrams \ref{under3}, \ref{under6}), 
  it needs to have 2 liberties and it will be able to
  capture the capturing chain, i.e.\
\item after replacing the cross-cutting stone the only remaining
  liberty of the capturing chain must not be a liberty of another
  chain so that the capturing chain can be saved by linking.
\end{itemize}
\[
\myblock{6\goheight}{
\board2--7/1--3/
\alphafalse
\whitetris{F1F2}
\whites{B1B2B3C3D3E2}
\blacks{C1C2D2D1E3F3G3G2G1}
\shipoutgoboard}{$\!\!\!\!\!$}{}{under1}
\hspace*{4\goheight}
\myblock{6\goheight}{
\board2--7/1--3/
\alphafalse
\blackfirstfalse
\whites{B1B2B3C3D3E2F2F1}
\blacks{C1C2D2D1E3F3G3G2G1}
\numbers{E1}
\shipoutgoboard}{$\!\!\!\!\!$}{}{under2}
\hspace*{4\goheight}
\myblock{6\goheight}{
\board2--7/1--3/
\alphafalse
\blackfirstfalse
\whites{B1B2B3C3D3E1E2F2F1}
\blacks{E3F3G3G2G1}
\numbers[2]{D2}
\shipoutgoboard}{$\!\!\!\!\!$}{}{under3}
\] \[
\myblock{6\goheight}{
\board2--8/1--3/
\alphafalse
\whitetris{G1G2}
\whites{B1B2B3C1C3D3E2F2}
\blacks{C2D2D1E1E3F3G3H3H2H1}
\shipoutgoboard}{$\!\!\!\!\!$}{}{under4}
\hspace*{3\goheight}
\myblock{6\goheight}{
\board2--8/1--3/
\alphafalse
\blackfirstfalse
\whites{B1B2B3C1C3D3E2F2G2G1}
\blacks{C2D2D1E1E3F3G3H3H2H1}
\numbers{F1}

\shipoutgoboard}{$\!\!\!\!\!$}{}{under5}
\hspace*{3\goheight}
\myblock{6\goheight}{
\board2--8/1--3/
\alphafalse
\blackfirstfalse
\whites{B1B2B3C1C3D3E2F1F2G2G1}
\blacks{E3F3G3H3H2H1}
\numbers[2]{D2}
\shipoutgoboard}{$\!\!\!\!\!$}{}{under6}
\]

\noindent
(If the stones \twtstone\ in diagrams \ref{under1}, \ref{under4} 
would all or individually be black then the positions would not be
semeai anymore but re-capture would still be possible.)

\end{document}